\documentclass[11pt]{amsart}

\usepackage{amssymb}
\usepackage{graphics}
\usepackage{graphicx}

\newtheorem{theorem}{Theorem}
\newtheorem{conj}{Conjecture}
\newtheorem{lemma}{Lemma}
\newtheorem{prop}{Proposition}

\newtheorem{cor}{Corollary}
\newtheorem{example}{Example}
\theoremstyle{definition}
\newtheorem{remark}{Remark}

\newcommand{\HF}{\widehat{HF}}
\newcommand{\on}{\operatorname}

\renewcommand{\d}{\partial}

\newcommand{\gr}{\on{gr}}

\newcommand{\x}{\tt{x}}
\newcommand{\p}{\tt{p}}
\newcommand{\q}{\tt{q}}
\renewcommand{\r}{\tt{r}}
\newcommand{\D}{\mathcal{D}}
\newcommand{\OO}{\mathcal{O}}
\newcommand{\uu}{{\bf u}}
\newcommand{\um}{{\bf u}_-}
\newcommand{\up}{{\bf u}_+}

\newcommand{\Z}{\mathbb{Z}}
\newcommand{\R}{\mathbb{R}}

\newcommand{\ps}{\ti{\psi}}
\newcommand{\ti}{\tilde}

\begin{document}

\author{Olga Plamenevskaya}
\address{Department of Mathematics, M.I.T., Cambridge, MA 02139}
\email{olga@math.mit.edu}
\title{Transverse knots and Khovanov homology}
%\subjclass{?}

\begin{abstract}
We define an invariant of transverse links in $(S^3,
\xi_{std})$ as a distinguished element of the Khovanov homology of
the link. The quantum grading of this invariant is the
self-linking number of the link. For knots, this gives a bound on
the self-linking number in terms of Rasmussen's invariant $s(K)$.
We prove that our invariant vanishes for transverse knot
stabilizations, and that it is non-zero for quasipositive braids.
We also discuss a connection to Heegaard Floer invariants.
\end{abstract}

\maketitle

\section{Introduction}
\subsection{Legendrian and Transverse knots}
There are two important classes of knots in a contact 3-sphere
$(S^3, \xi_{std})$: Legendrian knots and transverse knots.
Legendrian knots are everywhere tangent to the contact planes;
transverse knots are everywhere transverse to them. There are
simple "classical" invariants for both classes: the
Thurston--Bennequin and the rotation number for Legendrian knots,
and the self-linking number for transverse knots. While certain
knot types, e.g. all torus knots \cite{EH1}, are completely classified by
these invariants (in this case the knot type is said to be {\em
Legendrian resp. transversely simple}), for most knot types the
 classification is not
known. There exist smoothly isotopic Legendrian resp. transverse
knots with the same classical invariants which are not isotopic
through Legendrian resp. transverse knots. Legendrian knots are
somewhat better understood and enjoy a rich theory in the context
of symplectic field theory \cite{ElH}: a certain differential graded algebra
associated to a knot yields new Legendrian knot invariants.
Transverse non-simplicity of certain knot types was demonstrated
by Birman and Menasco in \cite{BM} and by Etnyre and Honda in
\cite{EH3}. However, the existing examples are sparse (those of
\cite{EH3} are not even explicit), and the proofs in \cite{BM} and
\cite{EH3} require a subtle analysis of braids and contact
manifolds. Unlike Legendrian knots, transverse knots do not have
any known efficient non-classical invariants.

\subsection{The invariant $\psi(L)$}
In this paper we introduce a transverse link invariant $\psi(L)$
as a distinguished element of  the Khovanov homology of $L$. Given
a closed braid diagram representing the transverse link, $\psi(L)$
is defined via a certain resolution of $L$. This invariant encodes
the self-linking number: the quantum degree of  $\psi(L)$ is given
by $sl(L)$. It also discerns transverse stabilizations: if $L$
arises as a transverse stabilization of another transverse link,
$\psi(L)$ vanishes. On the other hand, $\psi(L)\neq 0$ for
quasipositive braids. While we don't have any examples of
transverse knots distinguished by $\psi(L)$ but not $sl(L)$ (indeed, we show that
the invariant is the same for the pairs of transversely
non-isotopic knots from \cite{BM}), we hope that a connection to
Khovanov homology might be helpful. In particular, we establish
a bound on the self-linking number of a knot in terms of the Khovanov
homology knot invariant of Rasmussen \cite{Ra}.

\subsection{Khovanov homology and low-dimensional topology}
Khovanov homology is an invariant of knots and links introduced in
\cite{Kho}. Given a link $L$ in $S^3$, $Kh(L)$ is a graded
homology module whose graded Euler characteristic is the
unnormalized Jones polynomial of $L$. As was recently discovered,
the Khovanov homology has an interesting relation to
low-dimensional topology. Ozsv\'ath and Szab\'o \cite{OS}
construct a spectral sequence converging to the Heegaard Floer
homology $\widehat{HF}(Y)$ of the double cover $Y$ of $S^3$
branched over a link $L$; the $E^2$ term of this spectral sequence
is the (reduced) Khovanov homology of the link $L$. Rasmussen
\cite{Ra} uses Khovanov homology to give a combinatorial proof of
the Milnor conjecture (i.e., to determine the slice genus of a torus
knot).
 Our transverse link invariant suggests a
further connection to contact topology. In fact, we can define a
similar invariant in the reduced Khovanov theory; we conjecture
that this invariant is mapped  to the Ozsv\'ath-Szab\'o contact
invariant of the double cover of $(S^3,\xi_{std})$ branched over the
transverse link under the spectral sequence of \cite{OS}.

\subsection{Acknowledgements} I would like to thank Jake Rasmussen
for extremely helpful email correspondence, and Peter Kronheimer and
Ciprian Manolescu for illuminating conversations and encouragement.

\section{Preliminaries on transverse knots} It will be convenient
to work with closed braid representations of transverse knots.
Consider $S^3$ equipped with the (rotationally symmetric) standard
 contact structure  $\xi_{std}=\ker(dz-ydx+xdy)$.
It easy to see that any closed braid around $z$-axis can be made
transverse to the contact planes.  Moreover, by a theorem of
Bennequin \cite{Be} any transverse link in $(S^3, \xi_{std})$ is
transversely isotopic to a closed braid.

We adopt the usual notation for braid words. The braid group on
$b$ strings is generated by $\sigma_1, \dots, \sigma_{b-1}$, so
that $\sigma_i$ permutes the $i$-th and the $(i+1)$-th strings. We
will sometimes write a braid as a braid word, a certain product of
the generators $\sigma_1, \dots, \sigma_{b-1}$ and their inverses.
The positive resp. negative stabilization of a braid on $b$
strings is formed by adding the $(b+1)$-th string and multiplying
the braid word by $\sigma_b$ resp. $\sigma_b^{-1}$.

Of course, the same link can be represented by different braids.
The Markov theorem \cite{Bi} asserts that two braid words
describing the same link are related by a sequence of
stabilizations, destabilizations and conjugations in the braid
group (and, of course, the braid group identities). The Transverse
Markov Theorem  describes the relation between two braid
representations of the same transverse links.

\begin{theorem} \cite{W,OrSh} \label{markov} Let $L_1$, $L_2$ be two closed braids which
represent transversely isotopic links. Then $L_2$ can be obtained
from $L_1$ by a sequence of positive braid stabilizations and
braid isotopies.
\end{theorem}

The self-linking number $sl(L)$ is defined as follows. Fix a Seifert surface
$\Sigma$ for the link $L$, and let $v$ be a non-zero vector field in $T\Sigma\cap \xi$
along $L$ pointing out of $\Sigma$. Then, $sl(L)$ is the obsruction to extending
  $v$ over $\Sigma$.
Given a closed braid representing $L$, it can be computed as
\begin{equation}\label{sl-braid}
sl(L)=-b+n_+-n_-,
\end{equation}
where $b$ is the braid index, and $n_+$ and $n_-$ denote  the number of
positive resp. negative crossings.

The {\em stabilization} of a transverse link can be thought of as
negative braid stabilization. Unlike positive braid stabilization,
this operation changes the transverse type of the link: if $L'$ is
the result of stabilization of $L$, then
\begin{equation}
sl(L')=sl(L)-2.
\end{equation}

%For certain smooth knot types, the self-linking number completely
%classifies the transverse knots. If this is the case, the knot is
%called {\em transversely simple}.

%\begin{theorem}\label{trans_class}
%\end{theorem}

%Not all knot types are transversely simple, although the examples
%of non-isotopic transverse knots with the same linking number
%are sparse, and the proofs of non-simplicity are very subtle.
%Braid representations for ... are given by Birman and Menasco

%\begin{theorem}\label{birman-menasco}\cite{BM}
%The braids $K_1=\sigma_1^3\sigma_2...$ and  $K_2=\sigma_1^3\sigma_2^4...$
%have $sl(K_1)=sl(K_2)$, but are not isotopic as transverse knots.
%\end{theorem}

%Non-explicit examples are also provided by Etnyre and Honda \cite{EH}, who
%show that the $(2,3)$ cable of the right trefoil is not transversely
%simple: it admits two non-isotopic representatives with $sl=3$.

\section{Khovanov Homology}
In this section we give a brief review of Khovanov homology (more or less following
the review in \cite{Ra}). Unless otherwise specified, We
work with coefficients in $\Z$.

\subsection{Khovanov complex}
\begin{figure}[ht]
\includegraphics[scale=0.7]{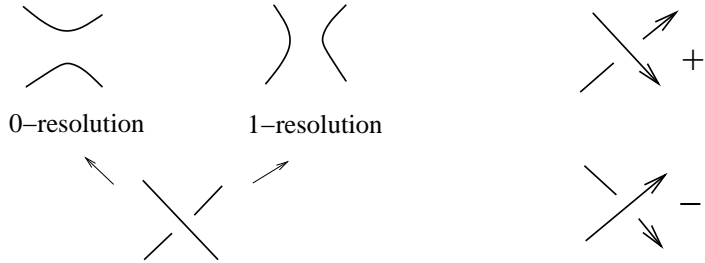}
\caption{Resolutions and signs of crossings.} \label{0-1-pos-neg}
\end{figure}

Given a link diagram $L$, we can resolve its crossings so that the
result is just the union of planar circles. Each crossing can be
resolved in two ways, called the 0-resolution and the 1-resolution
and shown in Fig.\ref{0-1-pos-neg}. Let $n$ be the number of
crossings of $L$; we will write $n=n_-+n_+$, where $n_+$ ($n_-$)
is the number of positive (negative) crossings. (See Fig.
\ref{0-1-pos-neg} for the usual sign conventions.)
 Then, complete
resolutions of $L$ can be conveniently labelled by  vertices of
the ``cube of resolutions'' $[0,1]^n$.
The underlying graded module for the Khovanov complex $CKh(L)$ is
the direct sum of $\Z$-modules associated to the vertices of
$[0,1]^n$,
$$
CKh(L)=\oplus_{v\in\{0,1 \}^n }CKh(L_v).
$$
Each  $CKh(L_v)$ is defined as follows. Let $U$ be the free graded
$\Z$-module generated by two elements, $\um$ and $\up$; the
grading $p$ is given by $p({\bf u}_\pm)=\pm 1$. Suppose that the
resolution $L_v$ consists of $k$ circles. We then set
$$
CKh(L_v)=U^{\otimes k}.
$$
In other words,  $CKh(L_v)$ is freely generated by $k$-tuples obtained by
labelling each circle in $L_v$ by either $\um$ or $\up$.

The module $CKh(L)$ is bi-graded.
The {\em homological} grading, which is constant  on each
 $CKh(L_v)$, is given by $\gr(v)=|v|-n_-$,
where $|v|$ is the number of 1's in the coordinates ov $v$.
Besides, there is the {\em quantum} grading $q(\uu)=p(\uu)+\gr(\uu)+n_+-n_-$.

Our next job is to describe the differential $d$ on $CKh(L)$. Loosely,
$d$ is the sum of maps $d_e$ associated to the edges of $[0,1]^n$.

Let $e$ be an edge of $[0,1]^n$, and denote by $v_e(0)$ resp.
$v_e(1)$ its initial resp. terminal end. The resolutions
$L_{v_e(0)}$ and $L_{v_e(1)}$ differ in one crossing only (which
is 0-resolved for $v_e(0)$ and 1-resolved for  $v_e(1)$), and
$L_{v_e(1)}$ is obtained from $L_{v_e(0)}$ in one of two ways:
either two circles of  $L_{v_e(0)}$ merge into one, or one circle
splits into two. In the first case, the map $d_e:CKh(L_{v_e(0)}
)\to CKh(L_{v_e(1)})$ is given by multiplication $m:U\otimes U\to
U$, where the two factors of $U\otimes U$ correspond to two
circles that merge, and the copy of $U$ in the image corresponds
to the resulting circle in $CKh(L_{v_e(1)})$. In the second case,
$d_e$ comes from  the comultiplication $\Delta:U\to U\otimes U$,
where the $U$ in the domain corresponds to the circle that splits.
It remains to define the maps $m$ and $\Delta$:
\begin{equation}
\begin{array}{ll}
m(\up\otimes\up)=\up & \Delta(\up)=\up\otimes\um+\um\otimes\up \\
m(\up\otimes\um)\otimes m(\um\otimes\up)=\um  & \Delta(\um)=\um\otimes\um.   \\
m(\um\otimes\um)=0
\end{array}
\end{equation}

Now, on the component  $CKh(L_v)$  the differential $d$ is defined by
$$
d=\sum_{e:v_e(0)=v}(-1)^{s(e)}d_e,
$$
where the sum is taken over all edges which have $v$ as their initial end.
The signs $(-1)^{s(e)}$ are chosen  so that $d^2=0$ (the choice is not unique,
but all the resulting chain complexes are isomorphic).

Khovanov \cite{Kho} shows that different diagrams for the same knot yield
 quasi-isomorphic chain complexes, so that the isomorphism classes of the
(bigraded) homology groups give an invariant of the link.
Actually, more is true: as conjectured in \cite{Kho} and proved in
\cite{Ja}, Khovanov's theory is functorial, and it follows that
there are honest homology groups, not just isomorphism classes,
associated to a link. We now turn attention to these the
functorial properties.

\subsection{Cobordisms and Invariance} \label{Kho-co} Given two links and an oriented cobordism
between them, there is an induced map between homology groups of the links. We briefly
describe this construction.

An oriented cobordism between two links  $L^0$ and $L^1$ is given
by an embedded smooth oriented compact surface $S$  in $\R^3\times
[0,1]$, such that $\d S=S\cap \d(\R^3\times[0,1])$, and $S\cap
(\R^3\times{i})=L^i$ for $i=0,1$. We may assume that $L^t= S\cap
(\R^3\times{t})$ is a link for all but finitely many values of $t$.
When $t$ passes through the critical value, the isotopy type of
the link changes by a {\em Morse move}, and the surface $S^t=
(S\cap\R^3\times[0,t])$ changes by an attachment of a handle (of index 0,
1, or 2). Further, we can fix a projection $\R^3 \to \R^2$, and
assume that it gives a regular projection for  $L^t$ for all but
finitely many special values of $t$ (and that the set of these
special values is disjoint from the set of the Morse critical
values). Thus, we obtain link diagrams (still denoted  $L^t$).
When $t$ passes through a special value where the projection of
the link is not regular, the link remains the same, but its
diagram changes by a Reidemeister move. The isomorphism class of the surface $S^t$
remains unchanged.

Therefore, the cobordism $S$ can be represented as a sequence of elementary cobordisms,
$$
S= S_1\cup S_2\cup\dots \cup S_k,
$$
where each cobordism $S_i$ between two diagrams $L^{t_i}$ and $L^{t_{i+1}}$ corresponds to either a Reidemester move
or a handle attachment.
Now, each $S_i$ induces a map $f_{S_i}:Kh(L^{t_i})\to Kh(L^{t_{i+1}})$. For Reidemeister moves,
$f_{S_i}$ comes from the quasi-isomorphisms between chain complexes  $CKh(L^{t_i})$
and $CKh(L^{t_{i+1}})$ mentioned in the previous section (we describe these quasi-isomorphisms in a little more
detail in Section \ref{def-psi}). For Morse moves,  $f_{S_i}$ is defined as follows \cite{Kho}.
We need two additional maps, $\iota:\Z\to U$ and $\epsilon:U\to \Z$, defined by
\begin{equation}
\begin{array}{ll}
\epsilon(\um)=1 & \iota(1)=\up \\
\epsilon(\up)=0
\end{array}
\end{equation}
Now, the attachment of a $0$-handle corresponds ro a ``birth'' of
a circle in the diagram, and the map on the chain complex is given
by $\iota$ (for all possible resolutions). Similarly, the
attachment of a $2$-handle (the ``death'' of a circle) gives the
map given by $\epsilon$. The attachment of 1-handle is given by
$m$ or $\Delta$ on each component of the chain complex, depending
on whether the 1-handle merges two circles of a particular
resolution or splits one circle into two. (Note that the
differential in Khovanov's theory is defined in a similar way: two
resolutions given by adjacent vertices of $[0,1]^n$ differ
precisely by the attachment of a 1-handle.) Finally, the map $f_S$
is defined as the composition of the maps induced by the
elementary cobordisms,
$$
f_S=f_{S_k}\circ\dots\circ f_{S_2}\circ f_{S_1}.
$$

Jacobsson  \cite{Ja} proves that up to a sign, the map $f_S$
depends on the isotopy class of $S\  rel\  \d S$ only, that is, if
$$
S= S_1\cup S_2\cup\dots \cup S_k \text{  and  } S= S'_1\cup S'_2\cup\dots \cup S'_{k'}
$$
are two decompositions of $S$ into elementary cobordisms, then
$$
f_{S_k}\circ\dots\circ f_{S_2}\circ f_{S_1}=\pm f_{S'_{k'}}\circ\dots\circ f_{S'_2}\circ f_{S'_1}.
$$
In particular, if two diagrams of a link are related by a sequence
of Reidemeister moves, then the induced isomorphism between the
homology groups is canonical up to a sign.

\section{Definition of the Invariant} \label{def-psi}

In this section we define the transverse link invariant $\psi(L)\in Kh(L)$.

First, we fix a braid diagram $L$ for our link, and
pick a distinguished element $\ps(L)$ in the chain complex $CKh(L)$.
We will check  that $\ps(L)$ is a cycle, so that it defines an
element $\psi(L)$ of the homology group $Kh(L)$. Finally, we show that
 $\psi(L)$ does not depend on the choice of the braid diagram
and remains the same under transverse link isotopies. This means that
  $\psi(L)$ is indeed an invariant of the transverse link.

Given a braid diagram  $L$ for our link,
  we choose a resolution  which is given by $b$ parallel strings:
that is, we take the 0-resolution for each positive crossing and
the 1-resolution for each negative crossing of $L$. Note that this is
the {\em oriented resolution} of the diagram; we denote it by
$L_o$.

 We set
\begin{equation} \label{defn}
\ps(L)=\um\otimes \um\otimes\dots \otimes \um\in U^{\otimes b}= CKh(L_o).
\end{equation}

\begin{prop} \label{cycle} The element $\ps(L)$ is a cycle in $(CKh(L), d)$.
\end{prop}

\begin{proof} The differential $d$ on $CKh(L_v)$ is the sum of maps
for all edges $e$ which have $v$ as their  initial end. By our choice of the
resolution $L_v$, such edges correspond to positive crossings. Moreover,
when a 0-resolution of a positive crossing is changed into a 1-resolution,
the two circles of $L_v$ which are ``connected'' by this crossing merge into one.
This means that each map $d_e$ is given by multiplication, and then
$$
d_e(\ps(L))=m(\um\otimes\um)=0.
$$
Taking the sum over all positive crossings, we see that $d(\ps(\D))=0$.
\end{proof}

\begin{prop} \label{deg} $\psi(L)\in Kh^{0,sl(L)}$.
\end{prop}

\begin{proof} By construction, $\psi(L)$ is a homogeneous element.
The homological and quantum gradings are easy to compute: since
the number of 1's for the chosen resolution is exactly the number
of negative crossings, $\gr(\psi(L))=0$. Now, $p(\psi(L))$ is the
braid index of $L$, and the formula $q(\psi(L))=sl(L)$ is an
immediate consequence of (\ref{sl-braid}) and the definition of
the quantum grading.
\end{proof}

\begin{figure} [ht]
\includegraphics[scale=0.7]{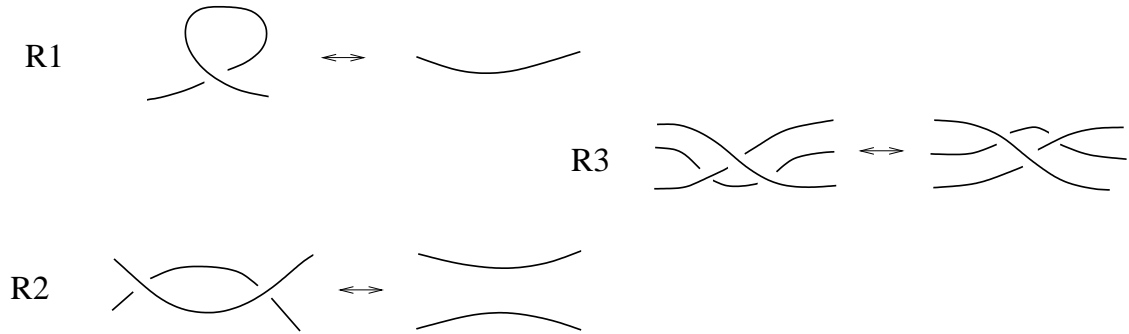}
\caption{Reidemeister moves in the transverse braid setting.}
\label{Rmoves}
\end{figure}

Now we want to check that $\psi(L)$ is independent of a particular braid representation of the transverse
knot. The Transverse Markov Theorem says that the braids representing two transversely isotopic knots
are related by a sequence
of positive stabilizations and braid isotopies. The two braid words will then be related
positive stabilizations, conjugations and the braid group identities. For our braid diagrams,
this yields a sequence of ``transverse Reidemeister moves'', as follows.
 Positive stabilization gives the move (R1) with
a positive crossing introduced (the other version of (R1) is
not allowed). The braid isotopies give the (R2) and (R3) moves,
all versions of which are allowed. In the following Lemma, we
check that the moves (R1)-(R3) respect $\psi(L)$ by analyzing the
effect of each move on Khovanov's homology. If  we were  dealing with
knots, it would suffice to consider only the versions of
(R1)-(R3) shown in Fig. \ref{Rmoves}, since all the other versions
can be obtained by a combination of these three. With braids, we
need to be more careful: there is another version of (R2) obtained
by turning our picture upside down; since we cannot turn braids
upside down, we actually need to consider both versions of (R2).
However, the two proofs are identical, so we only give one of
them. Also, it is not hard to check that all possible versions of
the (R3) move can be reduced to the one shown by a combination of
(R2) moves.

\begin{lemma} \label{R-invar} Let $L$ and $L'$ be two braid diagrams related by one of the three
transverse Reidemeister moves (R1), (R2), (R3), and denote by
$\rho_i:CKh(L)\to CKh(L')$, $i=1,2,3$, the associated quasi-isomorphisms between
the two chain complexes. Then
$$
\rho_i(\ps(L))=\pm \ps(L').
$$
\end{lemma}

\begin{proof} We  recall how the quasi-isomorphisms $\rho_i$ are constructed
in \cite{Kho}, and see what happens to the distinguished element $\ps(L)$.

(R1) move: The complex $CKh(L')$ decomposes as a direct sum $X_1\otimes X_2$,
where the  $X_2$ is acyclic, and $X_1$ is isomorphic to  $CKh(L)$ (Fig. \ref{R1-pf}). The isomorphism
$\rho_1$ is given by
$$
\begin{array}{ll}
\rho_1(\um)=\um\otimes\um \\
\rho_1(\up)=\up\otimes\um -\um\otimes\up,
\end{array}
$$
and we see that $\ps(L)$ is mapped to $\ps(L')$.

\begin{figure}[ht]
\includegraphics[scale=0.5]{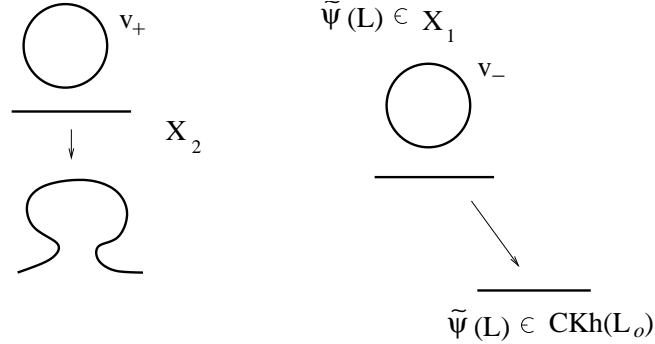}
\caption{Construction of the quasi-isomorphism $\rho_1$.}
\label{R1-pf}
\end{figure}
(R2) move: The four possible resolutions of the two extra crossings of $L'$
are shown in Fig. \ref{R2-pf}. Again , the complex  $CKh(L')$ decomposes as
 $X_1\otimes X_2\otimes X_3$, where $X_2$ and $X_3$ are both acyclic,
and $X_1$ isomorphic to  $CKh(L)$ via the isomorphism $\rho_2:CKh(L)\to CKh(L')$
given by
$$
\rho_2(x)=(-1)^{\gr(x)}(x+\iota(d_e(x))).
$$
(The map $\iota$ defined in Section \ref{Kho-co} and the map $\d_e$ corresponding to the edge $e$
 are shown
in the Figure, and the oriented resolution of $L$ is naturally identified with the oriented resolution of $L'$,
so that $x$ on the right-hand side actually lives in $CKh(L')$).

We see that  $\rho_2$ maps $\ps(L)$ to $\pm(\ps(L')+\iota(d_e(\ps(L')))$.
In the proof of Proposition \ref{cycle}, we've checked that
$d_e(\ps(L'))=0$. It follows that up to a sign, $\ps(L)$ is mapped to $\ps(L')$.

\begin{figure}[ht]
\includegraphics[scale=0.5]{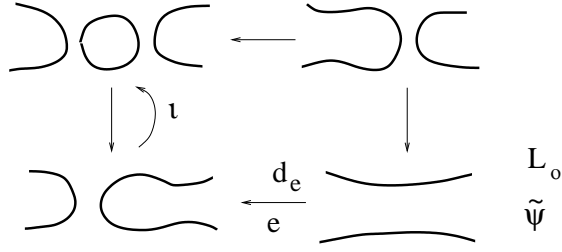}
\caption{Construction of the quasi-isomorphism $\rho_2$.}
\label{R2-pf}
\end{figure}

(R3) move: Now we have decompositions $CKh(L)= X_1\otimes
X_2\otimes X_3$ and $CKh(L')= X'_1\otimes X'_2\otimes X'_3$ where
$X_2$, $X_3$, $X'_2$, $X'_3$ are all acyclic, and there is an
isomorphism $\rho_3:X_1\to X'_1$.

We briefly describe how $X_1$ is formed. First, pick 1-resolutions
of crossings $\r$ and $\r'$ (note that the resulting diagrams are
isomorphic). With this fixed, consider all possible resolutions of
other crossings, and denote the direct sum of the associated
components of $CKh(L)$ resp. $CKh(L')$ by $CKh(L_{*1})$ resp.
$CKh(L'_{*1})$. (This is not a subcomplex.)
\begin{figure}[ht]
\includegraphics[scale=0.7]{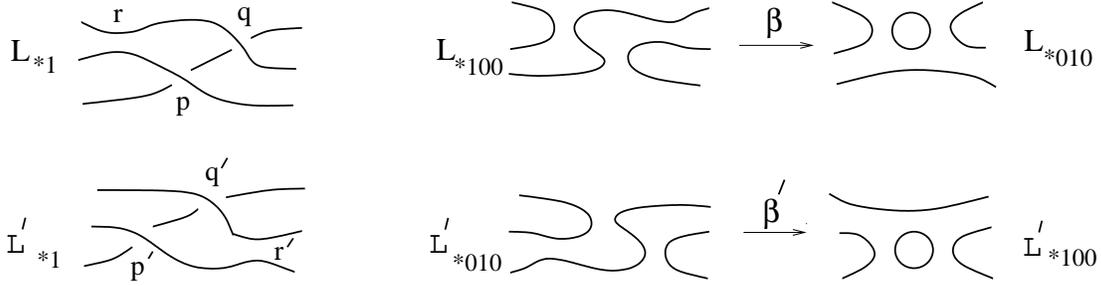}
\caption{Construction of complexes $X_1$ and $X'_1$.}
\label{R3-pf}
\end{figure}
 Next, denote by
$CKh(L_{*100})$ the part of $CKh(L)$ arising from all complete
resolutions of $L$ with a 1-resolution at $\p$ and 0-resolutions
at both $\q$ and $\r$; form $CKh(L_{*010})$, $CKh(L'_{*100})$ and
$CKh(L'_{*010})$ by analogy. Now, define

\begin{equation}
\begin{array}{l}
X_1=\{x+\beta(x)+y | x\in  CKh(L_{*100}), y \in CKh(L_{*1})\}\\
X'_1=\{x+\beta(x)+y | x\in  CKh(L'_{*010}), y \in CKh(L'_{*1})\},
\end{array}
\end{equation}
where $\beta:CKh(L_{*100})\to CKh(L_{*010})$ and
$\beta:CKh(L'_{*010})\to CKh(L'_{*100})$ are certain chain maps.
The isomorphism $\rho_3:X_1\to X'_1$ is given by
$$
\rho_3(x+\beta(x)+y)=x+\beta'(x)+y,
$$
where the natural identifications between $CKh(L_{*100})$ and
$CKh(L'_{*010})$ etc. are used. We do not describe the maps
$\beta$ and $\beta'$, referring the reader to \cite{Kho}: the only
thing we need to know is that $\ps(L)\in CKh(L_{*1})$ and $\ps(L')\in
CKh(L'_{*1})$, so $\rho_3(\ps(L))=\ps(L')$.

\end{proof}

We've checked that the distinguished element $\psi(L)\in Kh(L)$
behaves nicely under the three transverse Reidemeister moves, and
we know that any two transversely isotopic knots are related by a
sequence of such moves,  but why would an arbitrary transverse
isotopy between $L$ and $L'$ send $\psi(L)$ to $\psi(L')$? We have
to give the Transverse Markov Theorem another look: in \cite{W} it
is actually shown that an arbitrary transverse isotopy $S$ can be
smoothly modified  into a composition of braid isotopies and
positive stabilizations while the two links are fixed. Then, up to
an isotopy of $S\ rel\ \d S$, the cobordism $S$ between $L$ and
$L'$ decomposes as $S_1\cup\dots\cup S_k$, and Jacobsson's theorem
from Section \ref{Kho-co} implies that
$$
f_S(\psi(L))=f_{S_k}\circ\dots\circ f_{S_2}\circ
f_{S_1}(\psi(L))=\pm\psi(L').
$$
We have proved the following
\begin{theorem} The element $\psi(L)\in Kh(L)$ is an invariant of the
transverse link $L\in (S^3,\xi_{std})$, well-defined up to a sign.

\end{theorem}

\section{Properties of $\psi(L)$} \label{properties}

\subsection{Transverse stabilization}
\begin{theorem} \label{trst} If $L$ is the transverse stabilization of another transverse link,
then $\psi(L)=0$.
\end{theorem}

\begin{figure}[ht]
\includegraphics[scale=0.45]{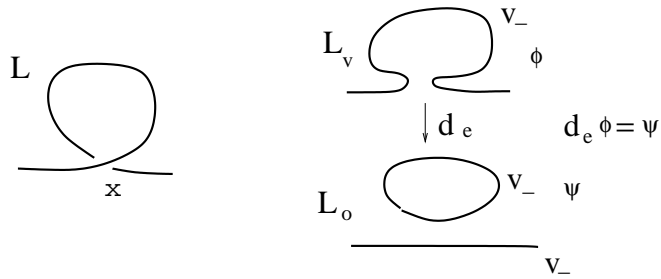}
\caption{Transverse stabilization and Khovanov's complex.}
\label{tra-sta}
\end{figure}

\begin{proof} We construct an element $\ti{\phi}\in CKh(L)$ such that $d \ti{\phi}=\ps(L)$.
Since $L$ is the result of a transverse stabilization (that is, an addition to the braid of an extra
string and an extra negative crossing $\x$),  it has a diagram with a
``negative kink''  as shown on Fig. \ref{tra-sta}. For the oriented resolution $L_o$, we take
the 1-resolution of the crossing $\x$. Let $e$ be the edge of the cube of resolutions corresponding to $\x$;
then $o$ is the terminal end of $e$. We denote by $v$ the initial end of $e$. In other words, we take
the 0-resolution of $\x$ to form $L_v$, and all the other crossings are resolved as in $L_o$.
Now, label all the components of $L_v$ by  $\um$, and set
$$
\ti{\phi}=\um\otimes\dots\um\in CKh(L_v).
$$
We compute $d \ti{\phi}$ as follows. The component $d_e:
CKh(L_v)\to CKh(L_o)$ of $d$ is given by comultiplication
$\Delta$, since the change of the resolution for $\x$ splits a
circle into two. Thus, $d_e \ti{\phi}=\ps$. Furthermore, similar
to proof of Proposition \ref{cycle},  all the other terms of $d$
on the component $CKh(L_v)$ correspond to positive crossings and
are given by multiplication maps, which send $\ti{\phi}$ to $0$.
It follows that $d \ti{\phi}=\ps$, as required.
\end{proof}

\subsection{Positive crossing resolution}
\begin{theorem}   \label{po-cro}
 Suppose that the transverse braid
$L^2$ is obtained  from the transverse braid $L^1$ by resolving a positive crossing
(note that it has to be a 0-resolution).
Let $S$ be the resolution cobordism, and $f_S:Kh(L^1)\to Kh(L^2)$ the associated map
on homology. Then
$$
f_S(\psi(L^1))=\pm\psi(L^2).
$$
\end{theorem}
\begin{proof} The cobordism $S$ is a composition of a 1-handle attachment and a Reidemeister move (R1),
as shown on Fig. \ref{pos-res}. On the component $CKh(L^1_o)$ of
the Khovanov's complex for $L^1$, the 1-handle attachment induces
a map given by comultiplication $\Delta$, since the handle splits
a circle on the oriented resolution. The map induced by the (R1)
move was analyzed in the proof of Lemma \ref{R-invar}. It follows
that the element $\ps(L^1)$, given by the $\um$ labels of all
circles for $L^1_o$, is mapped to $\ps(L^2)$ (given by the $\um$
labels on circles for $L^2_o$).

\begin{figure}[ht]
\includegraphics[scale=0.7]{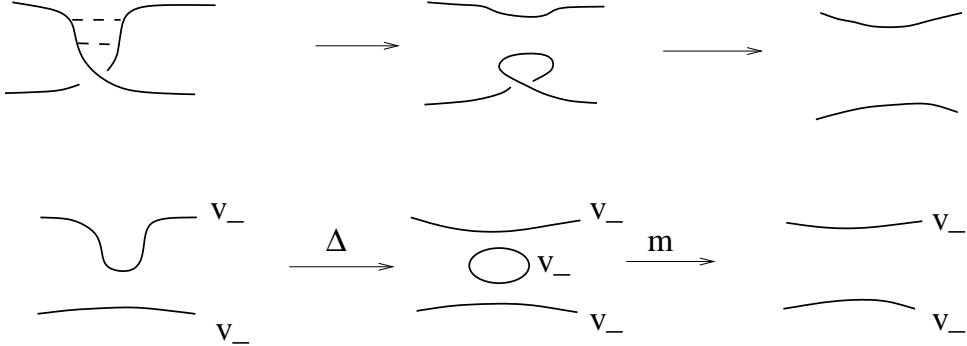}
\caption{Resolving a positive crossing.} \label{pos-res}
\end{figure}
\end{proof}

Recall that a braid is called {\em quasipositive} \cite{Ru} if its braid
word is a product of conjugates of the form $w \sigma_i w^{-1}$,
where $w$ is an arbitrary element of the braid group.

\begin{cor} \label{qua} If $L$ is represented by a quasipositive braid, then $\psi(L)\neq 0$. Moreover, it is
a primitive non-torsion element of $Kh(L)$.
\end{cor}

\begin{proof} Resolving a few positive crossings, we convert the braid representing $L$ into a braid
equivalent to a trivial one (of the same braid index). For the trivial braid $\OO$, there are no differentials
in the chain complex, and $\psi(\OO)$ is a generator of $Kh^{0,sl(\OO)}=\Z$. Since $\psi(L)$ is mapped to
$\psi(\OO)$, it must be non-torsion and primitive.
\end{proof}

\begin{remark} \label{pos-br} Let $L$ be a positive braid  of braid index $b$ with $n$ crossings.
Then, the homology of $L$ lies in non-negative homological degrees,
with  $Kh^{0, n-b}=\Z$,  $Kh^{0, n-b+2}=\Z$ and no other homology in $Kh^{0,*}$.
The element $\psi(L)$ is a generator of  $Kh^{0, n-b}$.
 \end{remark}

\begin{cor}\label{qua-non}  A transverse link $L$  represented by a quasipositive braid, then $\psi(L)\neq 0$ cannot be
obtained by a transverse stabilization of any other link.
\end{cor}

\begin{remark} Corollary  \ref{qua-non} follows easily from the fact that quasipositive
braid maximizes the self-linking number in its transverse link
type. More precisely, for an arbitrary transverse link $L$ the
slice-Bennequin inequality  \cite{Ru} gives
\begin{equation} \label{sl<g}
sl(L)\leq -\chi(\Sigma)
\end{equation}
where $\Sigma\subset B^4$ is a surface with boundary
 $\d \Sigma =L\in S^3=\d B^4$; for
quasipositive braids
 (\ref{sl<g}) becomes an equality. This bound was first proved
by Rudolph by means of gauge theory. It is interesting to note
that it can be obtained purely by Khovanov-homological methods.
Indeed, as in \cite{Ru}, it is straightforward to reduce the
question to the case of a positive braid representing a torus knot
$T_{p,q}$ (introducing
 positive crossings and
keeping track of how both sides of (\ref{sl<g}) change). Then, the
self-linking number is easily seen to be $2g(T_{p,q})-1$ (where
$g$ denotes genus), and $g_*(T_{p,q}=g(T_{p,q}))$. The last
identity is the Milnor conjecture, whose Khovanov homology proof
was recently obtained by Rasmussen \cite{Ra} (the original
gauge-theoretic proof is due to  Kronheimer and Mrowka).
 \end{remark}

%\begin{theorem} \label{c-sum} Let $L_1\#L_2$ be the connected sum
%of the transverse links $L_1$ and $L_2$. Then
%$p(\psi(L_1\#L_2))=\psi(L_1)\otimes\psi(L_2)$. Moreover, if
%$\psi(L_1)$ or $\psi(L_2)$ vanish, then $\psi(L_1\#L_2)$ also
%vanishes.

%\end{theorem}
\subsection{Negative crossings}
The following Proposition is useful for calculations and shows
that the invariant $\psi$ vanishes for many transverse links.

\begin{prop} \label{neg=0} Suppose that the transverse link $L$ is represented by a
closed braid whose braid word contains a factor of $\sigma_i^{-1}$ but no
$\sigma_i$'s for some $i=1, \dots, n$. (This means that all the crossings in the braid
diagram on the level between $(i-1)$-th and $i$-th string are
negative.) Then $\psi(L)=0$.
\end{prop}

\begin{proof} First of all, we delete all $\sigma_i^{-1}$ but one
from the braid word, obtaining a link that decomposes as a
connected sum of two links (connected by a negative crossing, the
$\sigma_i^{-1}$ that remains). 
\begin{figure}[ht]
\includegraphics[scale=0.7]{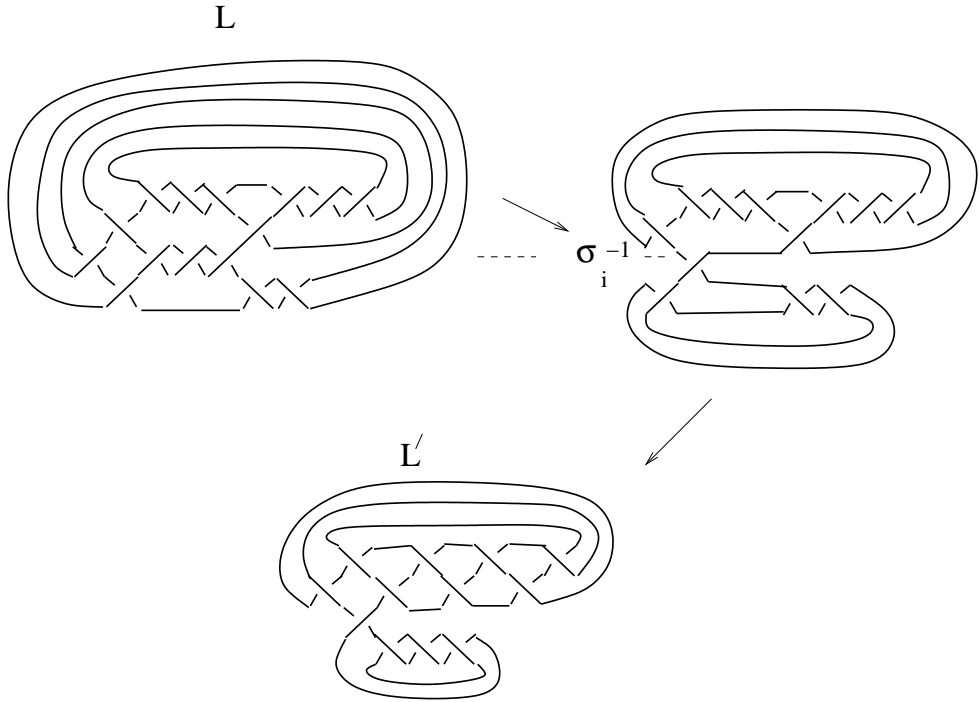}
\caption{From $L$ to $L'$.} \label{negcr}
\end{figure}
Then, we delete negative crossings
from and insert positive crossings into both components of the
connected sum, obtaining as a result two positive torus knots
connected by a negative crossing. This is illustrated on
Fig. \ref{negcr}. Denote the obtained  transverse
link by $L'$. By Theorem \ref{po-cro}, it suffices to show that
$\psi(L')=0$. Topologically, the link $L'$ is just the connected
sum of two torus knots, but its self-linking number is not maximal
(because we can connect the two components by a positive crossing
to increase $sl$). Connected sums of torus knots are transversely
simple \cite{EH2,Et}, so $L'$ is the transverse stabilization of
another link. By Theorem \ref{trst}, $\psi(L')=0$.
\end{proof}

\section{Examples}

\begin{example}
  For  $q>0$, $|p|\geq q$  let $L$
  be a transverse link
   of the $(p,q)$-torus
  link type.

\noindent (1) Suppose $p>0$. If $L$ maximizes the self-linking
number in its smooth link type, i.e. $sl(L)=pq-p-q$, then
$\psi(L)$ is a generator of $Kh^{0,pq-p-q}=\Z$. Otherwise
$\psi(L)=0$.

\noindent (2) If $p<0$,  then $\psi(L)$ vanishes.
\end{example}

\begin{proof} (1) We use transverse simplicity of positive torus links
\cite{EH1, Et}.
The (unique) transverse positive $(p, q)$-torus link with
$sl(L)=pq-p-q$ is represented by a positive braid with $q$ strings
and $p(q-1)$ crossings, so the result follows from Remark
\ref{pos-br}. If $sl(L)<pq-p-q$, then $L$ is obtained by
transverse stabilization, so $\psi(L)=0$ by Theorem \ref{trst}.

(2) Follows from Proposition \ref{neg=0}.
\end{proof}

\begin{example} Let the transverse links $L^1$ and $L^2$ be given
by the braids
$L^1=\sigma_1^{2p+1}\sigma_2^{2q}\sigma_1^{2r}\sigma_2^{-1}$ and
$L^2=\sigma_1^{2p+1}\sigma_2^{-1}\sigma_1^{2r}\sigma_2^{2q}$. 
It is shown in \cite{BM} that $L^1$ and $L^2$ are not
transversely isotopic when $p,q,r>1$, $p+1\neq q\neq r$ (although they are smoothly isotopic, and
$sl(L^1)=sl(L^2)$). However, we have
 $\psi(L^1)=\pm\psi(L^2)$. Indeed, both $\psi(L^1)$ and $\psi(L^2)$
 are generators of $Kh^{0,sl(L^i)}=\Z$.
\end{example}
\begin{proof} Because of  Proposition \ref{deg} and Corollary \ref{qua}, we
only need to check  that $Kh^{0,sl}=\Z$. The proof is
straightforward: the braids $L^1$ and $L^2$ have only one negative
crossing each, and for any link with this property,
the relevant part of $CKh(L)$ is easy to
understand. In the argument below,  $L$ will mean $L^1$ or $L^2$.
 For the given diagram $CKh^{-1, sl}(L)$ has
rank one and is generated by $\um \otimes \um$ in the complete resolution
of $L$ given by $0$-resolution of all crossings (this resolution
consists of two circles and is shown on Fig. \ref{0sl-pf}).
\begin{figure}[ht]
\includegraphics[scale=0.7]{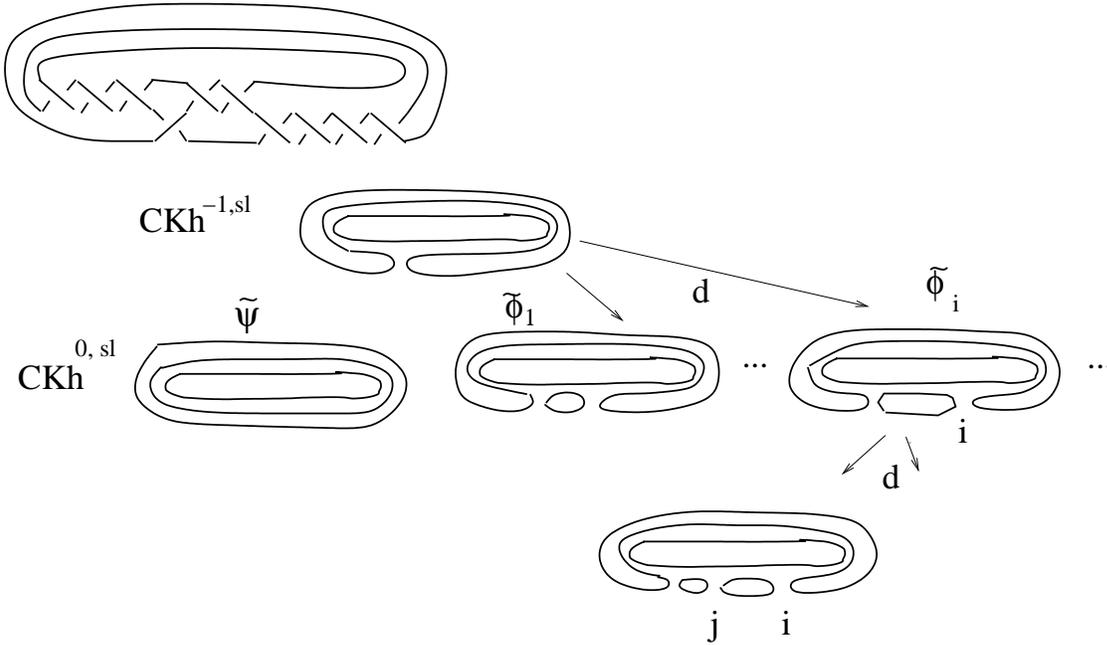}
\caption{Computing $Kh^{0,sl}(L)$. Each component of complete
resolutions shown must be labelled by  $\um$.} \label{0sl-pf}
\end{figure}
Similarly, $CKh^{0,sl}(L)$ has rank $2q$, with generators
$\ti{\phi}_i$ formed by the $\um$ labels on the complete
resolutions obtained as follows. The crossing that corresponds to
the $i$-th $\sigma_2$ in the product $\sigma_2^{2q}$ is
$1$-resolved, all other crossings are $0$-resolved. Further,
$d(\ti{\phi}_i)$  in turn comes from $1$-resolutions of the $j$-th
crossing in $\sigma_2^{2q}$ for all $j\neq i$. It follows that, apart
from $\ps(L)$, the only cycle in $CKh^{0,sl}(L)$ is
$\phi_1+\phi_2+\dots +\phi_{2q}$ (here we abuse our assumptions
and pretend that the coefficients are taken in $\Z/2\Z$; to be
honest, one needs to pick an appropriate choice of signs). The
latter cycle is the boundary of the generator of $CKh^{-1,
sl}(L)$, which finishes the proof.
\end{proof}

\section{A bound on the self-linking number}
In this section we obtain a bound on the self-linking number of a
transverse knot $K$ in terms of the knot invariant $s(K)$
introduced by Rasmussen \cite{Ra}. As we mention below, $s(K)$ is
defined as a certain quantum grading in Lee's version of the
Khovanov homology \cite{Lee}. (In Lee's construction, which works 
for rational coefficients only, the generators
for the complex $CKh'(K)$ and the gradings are the same as in
$CKh(K)$, but the differential is different.)
Rasmussen conjectures that this invariant 
coincides with the $s(K)$ invariant of Bar-Natan \cite{BN}, and is
twice the $\tau(K)$ invariant of Ozsv\'ath and Szab\'o
\cite{OStau}. Most importantly, $|s(K)|$ gives a lower bound of
the slice genus. If $K$ is  alternating, $s(K)$ is
simply the signature of the knot.

\begin{prop} \label{bound} For any transverse knot $K$
$$
sl(K)\leq s(K)-1.
$$
\end{prop}

\begin{proof}
As Jake Rasmussen pointed out to the author, the invariant $s(K)$
in \cite{Ra} is defined as $s(K)=\max q(\ti{x})+1$, where 
 $\ti x$ is an element of $CKh'(K)$ homologous to 
  the element that we have denoted $\ps$. (Note that Lee's differential 
does not preserve the quantum grading.)  Since $q(\ps)=sl(K)$
for a transverse knot $K$, it follows immediately  that 
$sl(K)\leq s(K)-1$.  
\end{proof}

\begin{remark} Proposition \ref{bound} gives an improvement for
the well-known Thurston--Bennequin \cite{Be} and slice--Bennequin
\cite{Ru} bounds on
$sl(K)$. The relation between Rasmussen's invariant and the 
slice--Bennequin inequality was independently established
by A. Shumakovich \cite{Sh}.
\end{remark}
\begin{cor} If $K$ is alternating, $sl(K)\leq \sigma(K)-1$, where
$\sigma(K)$ is the signature of the knot.
\end{cor}

\begin{remark} Since every bound for the self-linking number of
transverse knots is automatically a bound for $tb(K)+|r(K)|$, the
Thurston--Bennequin and rotation numbers of Legendrian knots
(Legendrian and transverse knots are related by push-offs
\cite{Et}), for an alternating Legendrian knot $K$ we have
$$
tb(K)+|r(K)|\leq \sigma -1.
$$
This bound was obtained in \cite{Pl1} via Heegaard Floer homology
techniques. Indeed, it is a special case of the inequality
$$
tb(K)+|r(K)|\leq 2\tau(K) -1,
$$
where $\tau(K)$ is the Ozsv\'ath--Szab\'o invariant \cite{OStau}.
The latter bound, together with Proposition \ref{bound},  supports
the conjecture about the  identity $s(K)=2\tau$, and gives yet
another connection between the Heegaard Floer theory and the
Khovanov homology.
\end{remark}

\section{Reduced homology and a relation to Ozsv\'ath--Szab\'o invariants}

A transverse link invariant can also be defined in the reduced
version of the Khovanov homology. We recall the construction of
the reduced complex. Starting with  a link $L$ with a marked point
on it, consider the usual Khovanov complex $CKh(L)$. For each
complete resolution of $L$, exactly one of the circles contains
the marked point. Let $CKh_{\um}(L)$ be the subcomplex generated
by those generators of $CKh(L)$ that have the label $\um$ on the
marked circle. Then, the reduced chain complex is the factor
$\widetilde{CKh}(L)= CKh(L)/CKh_{\um}(L)$, and $\widetilde{Kh}(L)$
is the corresponding homology. When the homology is taken with
coefficients in $\Z/2\Z$, $\widetilde{Kh}(L)$ is independent of
the choice of the marked point.

To define the transverse link invariant $\psi(L)$ in the reduced
homology, we take the same resolution as before, and pick a label
$\um$ for each unmarked circle and a $\up$ for the marked circle.
(Invariantly, we can consider the element
$$
\up\otimes\um\otimes \dots\otimes \um+
\um\otimes\up\otimes\dots\otimes\um+\dots+\um\otimes\um\otimes\dots\otimes\up
$$ in the non-reduced chain
complex, and then take its image in $\widetilde{CKh}(L)$.

The "reduced" invariant $\psi(L)\in\widetilde{Kh}(L)$ has the same
properties as those we proved in the non-reduced case (the proofs
would be identical).

In the Introduction, we mentioned the connection between
$\widetilde{Kh}(L)$ and the Heegaard Floer homology of the
3-manifold $\Sigma(L)$, which is a double cover of $S^3$ branched
over the link $L$ (here $L$ is a smooth link in $S^3$). More
precisely, there exist a spectral sequence whose $E^2$ term is
$\widetilde{Kh}(L)$, and $E^{\infty}$ term is the Heegaard Floer
homology $\HF(-\Sigma(L))$. When $L$ is an alternating link, this
spectral sequence collapses at the $E^2$ stage, so that
$\widetilde{Kh}(L)= \HF(-\Sigma(L))$.

When $S^3$ is equipped with the standard contact structure and $L$
is a transverse link,  the manifold $\Sigma(L)$ carries a natural
contact structure $\xi_L$ lifted from $S^3$.  In a related paper
\cite{Pl2}, we study the Ozsv\'ath--Sz\'abo contact invariant
$c(\xi_L)\in \HF(-\Sigma(L))$, associated to the contact
structure. It turns out that the properties of $c(\xi_L)$ are very
similar to those of $\psi(L)$; in particular, the results of
Section \ref{properties} hold true for $c(\xi_L)$. It is natural
to conjecture that the element $c(\xi_L)$ corresponds to $\psi(L)$
under the spectral sequence. For the case of alternating knots,
the conjecture is precise:

\begin{conj} Let $L$ be an alternating transverse link. Then $\psi(L)$
is mapped to $c(\xi_L)$ under the isomorphism
$\widetilde{Kh}(L)=\HF(-\Sigma(L))$.
\end{conj}

In \cite{Pl2}, we prove a special case of this conjecture.

\begin{theorem} Let $L$ be a transverse link represented by a
closed braid whose diagram is alternating. Then
$\psi(L)=c(\xi_L)$.
\end{theorem}

It should be said that the class of alternating braids is very
narrow, and that $\psi(L)=c(\xi_L)=0$ for most of them.

\end{document}